\documentclass[a4paper,11pt]{article}
\usepackage{graphicx}
\usepackage{amssymb}
\usepackage{latexsym, bm}
\usepackage{multicol}
\usepackage{indentfirst}
\usepackage{amssymb,amsfonts}
\usepackage{amsmath}
\newtheorem{theorem}{Theorem}
\newtheorem{lemma}{Lemma}

\title{New upper bound for multicolor Ramsey number of \\odd cycles
\thanks{Supported in part by NSFC(11671088), NSFFP(2016J01017) and CSC(201406655002).}}
\author{Qizhong Lin$^1$ and Weiji Chen$^2$\vspace*{0.3cm}\\
$^1$ Center for Discrete Mathematics, Fuzhou University \\Fuzhou 350108, China
\vspace*{0.3cm}\\
$^2$ College of Mathematics and Computer Science, Fuzhou University
\\Fuzhou 350108, China\vspace*{0.3cm}\\
 {\small\em Email: linqizhong@fzu.edu.cn}}
\date{}
\begin{document}
\maketitle
\begin{abstract}
Let $r_k(C_{2m+1})$ be the $k$-color Ramsey number of an odd cycle $C_{2m+1}$ of length $2m+1$.
It is shown that for each fixed $m\ge2$,
\[r_k(C_{2m+1})<c^{k}\sqrt{k!}\]
for all sufficiently large $k$, where $c=c(m)>0$ is a constant.
 This improves an old result by Bondy and Erd\H{o}s (Ramsey numbers for cycles in graphs, J. Combin. Theory Ser. B 14 (1973) 46-54).

\medskip

{\bf Keywords:} \  Ramsey number; odd cycle; upper bound
\end{abstract}

\section{Introduction}

Let $G$ be a graph. The multicolor {\em
Ramsey number} $r_k(G)$ is defined as the minimum integer $N$ such
that each  edge coloring of the complete graph $K_N$ with $k$ colors contains a monochromatic $G$ as a subgraph.
The Tur\'{a}n number $ex(N;G)$ is the maximum number of edges among all graphs of order
$N$ that contain no $G$. For the complete bipartite graph $K_{t,s}$ with $s\ge t$,
a well known argument of K\"ov\'ari,  S\'os, and Tur\'an \cite{kst-54} gives that
$
ex(N;K_{t,s})\le \frac{1}{2}\left[(s-1)^{1/t}N^{2-1/t}+(t-1)N\right].
$ For large $N$, the upper bound was improved by F\"uredi \cite{fur}
to $\frac{1}{2}((s-t+1)^{1/t}+o(1))N^{2-1/t}$.
Let $N=r_k(K_{t,s})-1$. Since there exists a $k$-coloring of the edges of $K_N$ such that it contains no monochromatic $K_{t,s}$,
which implies that each color class can have at most $ex(N;K_{t,s})$ edges. Thus ${N\choose 2}\le k\cdot ex(N;K_{t,s})$.  From an easy calculation, we have $r_k(K_{t,s})\le(s-t+1+o(1))k^t$ as $k\to\infty$. Hence $r_k(G)$ can be bounded from above by a
polynomial of $k$ if $G$ is a bipartite graph.

However, the situation becomes dramatically different when $G$ is non-bipartite. Denote $r_k(K_3)$ by $r_k(3)$ for short.
An old problem proposed by Erd\H{o}s is to determine $$\lim_{k\to\infty}(r_k(3))^{1/k}.$$
It is known from Chung \cite{chung} that $r_k(3)$ is super-multiplicative in $k$ so that $\lim_{k\to\infty}(r_k(3))^{1/k}$ exists.
Up to now, we only know that
\[
1073^{k/6}\le r_k(3)\le c\cdot k!,
\]
where $c>0$ is a constant, see \cite{a-h,Chung-g,Exoo,Fr-s,wan} and their references for more details.

Let $C_{2m+1}$ be an odd cycle of length $2m+1$.
For $m=1$, the multicolor {Ramsey number} $r_k(3)$ has attracted a lot of attention.
For general fixed integer $m\ge2$, Erd\H{o}s and Graham \cite{e-g} showed that
\begin{eqnarray}\label{b-e}
m2^k< r_k(C_{2m+1})< 2(k+2)!m.
\end{eqnarray}
Bondy and Erd\H{o}s \cite{b-e} observed that
\begin{eqnarray}\label{b-e}
m2^k+1\le r_k(C_{2m+1})\le(2m+1)\cdot(k+2)!.
\end{eqnarray}
For the lower bound, a recent result by Day and Johnson \cite{d-j} gives that for $m\ge2$, there exists a constant $\epsilon=\epsilon(m)>0$ such that $r_k(C_{2m+1})> 2m\cdot(2+\epsilon)^{k-1}$ for all large $k$.
For the upper bound, which was improved by Graham, Rothschild and Spencer
\cite{graham-roth-spen} to $r_k(C_{2m+1})<2m\cdot(k+2)!$.
In particular, for $m=2$, Li \cite{li-jgt} showed that $r_k(C_5)\le c\sqrt{18^kk!}$  for all
$k\ge3$, where $0<c<1/10$ is a constant. However, there are not too many substantial progress of $r_k(C_{2m+1})$ for $m\ge3$.

Let us point out that the situation is much different when $k$ is fixed.
For $k=2$, Bondy and Erd\H{o}s \cite{b-e}, Faudree and Schelp \cite{f-s(Cycle-cycle)} and Rosta \cite{ro} independently obtained that $r_2(C_{2m+1})=4m+1$ for all $m\ge2$. For $k=3$, {\L}uczak [9] proved that $r_3(C_{2m+1})=(8+o(1))m$ as $m\to\infty$ by using the regularity lemma.
Kohayakawa, Simonovits and Skokan \cite{k-s-s} used {\L}uczak¡¯s method together with stability
methods proved that $r_3(C_{2m+1})=8m+1$ for sufficiently large $m$.
Recently, Jenssen and Skokan \cite{j-s} established that $r_k(C_{2m+1})=2^km+1$ for all fixed $k$ and
 sufficiently large $m$.

In this short note, we have an upper bound for $r_k(C_{2m+1})$ as follows.
\begin{theorem}\label{odd-cycle}
Let $m\ge2$ be a fixed integer. We have
\[r_k(C_{2m+1})<c^{k}\sqrt{k!}\]
for all sufficiently large $k$, where $c=c(m)>0$ is a constant.
\end{theorem}

{\em Remark.} We do not attempt to optimize the constant $c=c(m)$ in the above theorem, since we care more about the
exponent of $k!$.

Let $N=r_k(G)-1$. From the definition, there exists a $k$-edge coloring of $K_N$ containing no monochromatic $G$.
In such an edge coloring, any graph induced by a monochromatic set of edges is called a Ramsey graph.
Let $\epsilon>0$ be a constant.
Under the assumption that each Ramsey graph $H$ for $r_k(C_{2m+1})$ has minimum degree at least $\epsilon d(H)$ for large $k$,
 Li \cite{li-jgt} showed that $r_k(C_{2m+1})\le \big(c^kk!\big)^{1/m}$,
where $d(H)$ is the average degree of $H$ and $c=c(\epsilon,m)>0$ is a constant.

\section{Proof of the main result }

In order to prove Theorem \ref{odd-cycle}, we need the following well-known result.

\begin{theorem}(Chv\'atal \cite{chv}) \label{chv}
Let $T_m$ be a tree of order $m$.  We have
\[
r(T_m,K_n)=(m-1)(n-1)+1.
\]
\end{theorem}

For a graph $G$, let $\alpha(G)$ denote the independence number of $G$.
\begin{lemma}(Li and Zang \cite{l-z})\label{l-z}
Let $m\ge2$ be an integer and let $G=(V,E)$ be a graph of order $N$
that contains no $C_{2m+1}$. We have
\[\alpha(G)\ge \frac{1}{(2m-1)2^{(m-1)/m}}\Bigg(\sum_{v\in V}d(v)^{1/(m-1)}\Bigg)^{(m-1)/m},\]
where $d(v)$ is the degree of $v$ in graph $G$.
\end{lemma}

\medskip
\noindent{\bf Proof of Theorem \ref{odd-cycle}.} Let $m\ge2$ and
$k\ge3$ be integers. For convenience, let $r_k=r_k(C_{2m+1})$ and $N=r_k-1$. Let $K_N=(V,E)$ be the complete
graph on vertex set $V$ of order $N$. From the definition, there exists an
edge-coloring of $K_N$ using $k$ colors such that it contains no
monochromatic $C_{2m+1}$. Let $E_i$ denote the monochromatic set of edges in color $i$ for $i=1,2,\dots,k$. Without loss of generality,
we may assume that $E_1$ has the largest cardinality among all
$E_i\,^{'}s$. Therefore $|E_1|\ge{N\choose 2}\big/k$. Let $G$
be the graph with vertex set $V$ and edge set $E_1$. Then the
average degree $d$ of $G$ satisfies
\[d=\frac{2|E_1|}{N}\ge\frac{N-1}{k}=\frac{r_k-2}{k}.\]

Consider an independent set $I$ of $G$ with $|I|=\alpha(G)$. Since
any edge of $K_N$ between two vertices in $I$ is colored by one of
the colors $2,3,\dots,k$, the subgraph induced by $I$ is an
edge-colored complete graph using $k-1$ colors, which contains no
monochromatic $C_{2m+1}$. Thus $|I|\le r_{k-1}-1$, and thus Lemma \ref{l-z} implies that
\begin{equation}\label{0}
r_{k-1}-1 \ge a \left(\sum_{v\in V} d(v)^{1/(m-1)} \right)^{(m-1)/m},
\end{equation}
where $a=a(m)$ is a constant.

\medskip
\noindent{\bf Claim.}  We have that
\begin{equation*}\label{1}
r_k\le c_1\sqrt{k}\,\, r_{k-1}
\end{equation*}
for some constant $c_1=c_1(m)$.

\noindent{\bf Proof.} For $m=2$, the assertion is clear since the inequality (\ref{0}) implies that
$$r_{k-1}-1 \ge a \sqrt{Nd }\ge a \sqrt{\frac{(r_k-1)(r_k-2)}{k}}> a \frac{r_k-2}{\sqrt{k}}.$$
In the following, we shall suppose $m\ge3$ and separate the proof into two cases.

{\bf Case 1.} The maximum degree $\Delta(G)$ of the graph $G$ satisfies $\Delta(G)>\frac{r_k}{\sqrt{k}}$, i.e.
there is some vertex $v$ such that $d(v)>  \frac{r_k}{\sqrt{k}}$. As the neighborhood $N(v)$ of $v$ contains no path $P_{2m}$ of order $2m$, we have from Theorem \ref{chv} that
\[
r_{k-1} -1 \ge \alpha(G)\ge  \frac{d(v)}{2m} > \frac{r_k}{2m\sqrt{k}},
\]
and so the claim holds for Case 1.

{\bf Case 2.} $\Delta(G)\le \frac{r_k}{\sqrt{k}}$.
Define a function
\[
f(x_1,x_2,\dots,x_N)=\sum_{i=1}^N x_i^{1/(m-1)},
\]
and consider the following optimization problem
\[
\left\{ \begin{array}{c}
 \min f =\min f(x_1,\dots,x_N), \\
   s.t. \;\; \sum_{i=1}^N x_i =Nd, \\
  \text{and}\;\; 0\le x_i \le \frac{r_k}{\sqrt{k}}\;\;\; \text{for}\;\; 1\le i\le N.
\end{array} \right.
\]
Using Lagrange multiplier method by setting 
$$L=L(x_1,\dots,x_N,\lambda)=f(x_1,\dots,x_N)+\lambda\Big(\sum_{i=1}^N x_i -Nd\Big),$$
 we find the {\em unique extreme point} ${\bf x}=(d,\dots,d)$.
Note that the Hessian matrix of $L$ (also $f$) at the point ${\bf x}$ is negative definite since its diagonal elements equal
$\frac{2-m}{(m-1)^2}d^{(3-2m)/(m-1)}$ which is negative for $m\ge3$ while the off diagonal elements equal zero,
so $f$ takes the maximum value at ${\bf x}$. However, the point ${\bf x}$ is not what we want.
Let
$$D=\Big\{(x_1,\dots,x_N): \; 0\le x_i \le \frac{r_k}{\sqrt{k}},\; 1\le i\le N\Big\}$$
 denote the feasible region of the above optimization problem. 
 
 Note that $f$ is a concave and continuous function with $D$ closed,
 hence the point we shall find such that $f(x_1,\dots,x_N)=\min f$ must be at the boundary of $D$, namely, at least
one $x_i=0$ or  $x_i=\frac{r_k}{\sqrt{k}}$, say $x_N=0$ (The case that $x_i=\frac{r_k}{\sqrt{k}}$ is similar). Hence the optimization problem become that for $N-1$ variables $x_1,\dots,x_{N-1}$.
By induction, we see that $f$ attains the minimum value at the point which has as many $x_i=0$ (or $x_i=\frac{r_k}{\sqrt{k}}$) as possible.
Let $$h=\Big\lfloor\frac{Nd\sqrt{k}}{r_k}\Big\rfloor.$$
Thus, 
we may take $x_i=\frac{r_k}{\sqrt{k}}$ for $1\le i\le h$, $x_i=0$ for $h+2\le i\le N$ and $x_{h+1}=Nd-h\frac{r_k}{\sqrt{k}}$,
and $f$ attain the minimum value at $(\frac{r_k}{\sqrt{k}},\dots,\frac{r_k}{\sqrt{k}},x_{h+1},0,\dots,0)$.
That is to say, 
\begin{align}\label{minf}
\nonumber\min f &=f\Big(\frac{r_k}{\sqrt{k}},\dots,\frac{r_k}{\sqrt{k}},x_{h+1},0,\dots,0\Big)
\ge h \Big( \frac{r_k}{\sqrt{k}}\Big)^{1/(m-1)} \\ &= 
\Big\lfloor\frac{Nd\sqrt{k}}{r_k}\Big\rfloor \Big( \frac{r_k}{\sqrt{k}}\Big)^{1/(m-1)}
 \ge\frac{1}{2}Nd \Big( \frac{r_k}{\sqrt{k}}\Big)^{-(m-2)/(m-1)}.
\end{align}
Therefore, from (\ref{0}) and (\ref{minf}), we obtain
\[
r_{k-1}-1 \ge \frac{a}{2} \left[ Nd \Big( \frac{r_k}{\sqrt{k}}\Big)^{-(m-2)/(m-1)} \right]^{(m-1)/m}.
\]
As $Nd \ge (r_k-1)\frac{r_k-2}{k}> \frac{(r_k-2)^2}{k}$, we have
\[
r_{k-1}-1 > \frac{a}{2} \left[ \frac{(r_k-2)^2}{k} \Big( \frac{r_k}{\sqrt{k}}\Big)^{-(m-2)/(m-1)} \right]^{(m-1)/m} 
> \frac{a}{2}\cdot \frac{r_k-2}{\sqrt{k}}.
\]
This completes the proof of Case 2 and hence the claim.\hfill$\Box$

Note that $r_2(C_{2m+1})=4m+1$ for $m\ge2$, see \cite{b-e,f-s(Cycle-cycle),ro}, and
repeatedly apply the above claim yields that
\begin{align*}r_k&\le c_1 \sqrt{k}\,\, r_{k-1}\le
c_1^{k-2}\sqrt{k(k-1)\cdots3}\cdot r_2(C_{2m+1})<c^k\sqrt{k!}
\end{align*}
for some constant $c=c(m)$.
This completes the proof of Theorem \ref{odd-cycle}. \hfill$\Box$

\medskip

{\bf Acknowledgment.}  We are grateful to the referees for giving detailed and very invaluable suggestions and comments
that improve the presentation of the manuscript greatly.

\end{document}